\documentclass[10pt]{amsart}
\usepackage[a4paper,top=2cm,bottom=2.5cm,left=3.5cm,right=3.5cm]{geometry}
\usepackage{enumerate} 
\usepackage{booktabs} 
\usepackage[utf8]{inputenc} 
\usepackage[T1]{fontenc} 
\usepackage{avant} 
\usepackage{times} 
\usepackage{amssymb,amsthm}
\usepackage{mathrsfs}
\usepackage[english]{babel} 
\usepackage{lmodern}
\usepackage{textcomp}
\usepackage{bm}
\usepackage{mathtools}
\usepackage{xfrac}
\usepackage{graphicx}
\usepackage{xcolor}
\usepackage[font=footnotesize,tableposition=top,figureposition=bottom,skip=0pt]{caption} 
\usepackage{hyperref}
\hypersetup{%
	linkcolor=blue,
	anchorcolor=black,
	citecolor=red,
	filecolor=magenta,
	menucolor=red,
	urlcolor=red,
	colorlinks=false,
	bookmarksnumbered=true,
	bookmarksopen=true,
	bookmarksopenlevel=0,
	breaklinks=true
}
\newtheorem{thm}{Theorem}

\newtheorem{cor}{Corollary}

\theoremstyle{definition}

\theoremstyle{remark}

\numberwithin{equation}{section}

\begin{document}
	
	\title[OSCILLATING SEQUENCES and QUASI-DISCRETE SPECTRUM]{OSCILLATING SEQUENCES OF HIGHER ORDERS AND TOPOLOGICAL
		SYSTEMS OF QUASI-DISCRETE SPECTRUM
		}
	\author{AIHUA FAN }
	\address{Fan Ai-Hua: LAMFA UMR 7352, CNRS, University of Picardie, 33, rue Saint Leu, 80039 Amiens CEDEX
		1, France}
	\email{\href{mailto:ai-hua.fan@u-picardie.fr}{ai-hua.fan@u-picardie.fr}}
	\thanks{}
		
	\begin{abstract}
		Fully oscillating sequences are orthogonal to all topological dynamical systems of quasi-discrete spectrum
		in the sense of Hahn-Parry. This orthogonality concerns with not only simple but also multiple ergodic means. 
		It is stronger than that required by Sarnak's conjecture .
	\end{abstract} 
	\maketitle
	

	
	
	\section{Main Statements}
	In this note  
	\footnote{
		This note was followed by Ai Hua FAN, Fully oscillating sequences and weighted multiple ergodic limit,
		C. R. Acad. Sci. Paris, Ser. I355 (2017), 866-870. Using different methods, both notes give almost the same result. As asked by some authors, I am putting this note on ArXiv.  
	}, 
	 we introduce the notion of  oscillating sequences of higher orders 
	and point out its relation to topological dynamical systems of quasi-discrete spectrum 
	in the sense of Hahn-Parry \cite{14hahn1965minimal}. This relationship is an orthogonality between oscillating 
	sequences of all orders and systems of quasi-discrete spectrum. This orthogonality is stronger than that required by
	Sarnak's conjecture (see \cite{24Sarnak2009three,25sarnak2010mobius} for Sarnak's conjecture). 
	Actually it is proved that {\em  multiple ergodic averages}  along polynomial times and weighted by a 
	fully oscillating sequence tend to zero everywhere. The oscillation may be reduced to some  finite order
	if the system has a quasi-discrete spectrum of finite order.
	
	A
     sequence of complex numbers $ c = (c_n)_{n\ge 0} $ is said to be \textit{oscillating of order $ d $} ($ d\ge 1 $) 
	 if for any real polynomial $ P\in \mathbb{R}_d[z] $  of degree less than or equal to $ d $ we have
	 \begin{equation}\label{eq:OSCd}
	 \lim\limits_{N\to \infty}\frac{1}{N}\sum_{n=0}^{N-1}c_{n}e^{2\pi i{P}(n)} = 0.
	 \end{equation}
	Then we write $ (c_n) \in \mathrm{OSC}_d $. The set $ \mathrm{OSC}_d $ is decreasing as $d$ increases to infinity. If
	$ (c_n) $ is oscillating of all orders $ d\ge 1 $, we say that it is\textit{ fully oscillating}. Then 
	we write $ (c_n) \in \mathrm{OSC}_{\infty} $. We say that $ (c_n) $ is \textit{oscillating of exact order $ d $} if it 
	is oscillating of order $ d $ but not oscillating of order $ d + 1 $. Then we write $ (c_n) \in \mathrm{OSC}^*_d $.
	Recall that the oscillation of order one means nothing but that the Fourier-Bohr spectrum of the sequence is empty. 
	Therefore the oscillations of higher orders are improvements of the fact that the Fourier-Bohr spectrum
	is empty. The notion of oscillation of order $1$ appeared in \cite{10fan2015oscillating}.
	
    Let us make a precise statement.  In the following, $ d^{\circ}q $ denotes the 
	degree of a polynomial $ q $.
	\begin{thm}\label{thm:1}
		Let $ \ell\ge 1 $ be an integer. Assume that
		
	\indent \textnormal{(i)}\, $ q_1,\cdots,q_{\ell}\in \mathbb{Q}[z] $ are $ \ell $ polynomials taking values in $ 
	\mathbb{N} $  with $ \Delta := \max(d^\circ q_1, \cdots, d^\circ q_\ell) $;
	 
	 \indent \textnormal{(ii)}  $ (X,T) $ is a dynamical system of quasi-discrete spectrum of order $ d_T\,(1\le d_T 
	 \le \infty) $;
	 \indent \textnormal{(iii)} $ (c_n) $ is an oscillating sequence of order $ d_T\Delta $ such that $ 
	 \sum_{n=1}^{N}|c_n|= O(N) $.\newline
	Then for any continuous function $ F \in C(X^\ell) $ and any point $ x\in X $, we have
	\begin{equation}\label{eq:1-1}
	\lim\limits_{N\to\infty}\dfrac{1}{N}\sum_{n=0}^{N-1}c_n F(T^{q_1(n)}x, \cdots, T^{q_\ell(n)}x) = 0.
	\end{equation}
	\end{thm}

	\medskip
	

     We need recall the notion of quasi-discrete spectrum, due to Hahn-Parry \cite{14hahn1965minimal} which is based on Abramov's theory on measure-theoretic 
    dynamics \cite{1abramov1962metric} using the concept of quasi-eigenfunction due to Halmos and von Neumann \cite{15halmos1942operator}. 
    In the present note, not as in \cite{14hahn1965minimal}, we don't assume that our dynamical systems are totally minimal, but transitive. 
    
     By a dynamical system we mean a pair $ (X,T) $ where $ X $ is a compact metric space and
    $ T\colon X\to X $ is a continuous map. Assume that $ (X,T) $ is transitive, i.e. 
    the orbit $ O(x) := \{ T^n x : n \ge 0 \}$ of some $ x \in X $ is dense in X. Let $ C(X)$ be the 
    Banach algebra of continuous complex valued functions on $ X $ and let
    \[
    G(X) := \{f \in C(X) : \forall x \in X,  |f(x)| = 1 \} 
    \]
   which is  a group under multiplication. All groups appearing below are subgroups of $G(X)$. 
    The transformation $ T $ induces 
    an isometry $ V_T $ of the algebra $ C(X) $, defined by $ V_T f(x) = f(Tx)$. 
    We say that $ f \in C(X)$, $f \ne 0$, is an \textit{eigenfunction} if there is a complex number $\lambda \in \mathbb{C} $ 
    for which
    \begin{equation}\label{eq:1-4}
    f\circ T= \lambda f.
    \end{equation}
    The number $ \lambda $ is called an \textit{eigenvalue}. Let $ H_1 $ be the group of all eigenvalues. 
    The eigenfunctions corresponding to the eigenvalue $1$ are called {\em  invariant functions}. The transitivity of 
    $T$ implies that invariant functions are constant functions, and $H_1 \subset K$ where $ K $ is the group
    $\{z \in \mathbb{C} : |z| = 1\}$ under multiplication, and eigenfunctions have constant modulus. 
    We define $G_1$ to
    be the group of all eigenfunctions $f$ belonging to $G(X)$. 
    We see that $H_1 \subset G_0 \subset G_1$, where we make the convention 
     $G_0: = H_1$.
    
    Now we define $H_n$ and $G_n$ ($n \ge 2$) by induction. Assume $H_1 \subset H_2 \subset \cdots \subset H_n $ and
    $ G_0\subset G_1 \subset G_2 \subset \cdots \subset G_n$ have already been defined such that $ H_{i+1} \subset G_i$ for 
    $0\le  i < n $. 
    Let $G_{n+1} $ be the set of all $ f_{n+1} \in G(X)$ 
    such that there is a $ g_n \in  G_n $ with 
    \begin{equation}\label{eq:1-5}
    f_{n+1}\circ T = g_n f_{n+1}.
    \end{equation}
Let $H_{n+1}$ be the set of all $g_n \in G_n $ for which there is an $f_{n+1} \in  G_{n+1}$ verifying   
\eqref{eq:1-5}. Denote
\[
G := \bigcup_{n=1}^{\infty}G_{n}, \quad H := \bigcup_{n=1}^{\infty}H_{n}.
\]
The elements of $H$ are called \textit{quasi-eigenvalues} and those of $G$ are called \textit{quasi-eigenfunctions}.

A dynamical system $(X, T)$ is said to have \textit{quasi-discrete spectrum} if the algebra generated
by the quasi-eigenfunctions is dense in $C(X)$, or equivalently the linear span of quasi-eigenvalues
is dense in $C(X)$ because $G$ is a multiplicative group. By using Stone-Weierstrass
theorem, this is equivalent to saying that quasi-eigenfunctions separate points of
$ X $. If, furthermore, $ G_d = G_{d+1} $ and $d_T$ is the least such integer, we say that $(X,T)$ has \textit{quasi-discrete
spectrum of order} $d_T$.

\medskip     
It was proved in \cite{17hoare1966affine2} that all minimal affine transformations of any connected compact
abelian group have quasi-discrete spectrum. This, together with Theorem~\ref{thm:1}, implies immediately
the following result, the conclusion of which is stronger than 
Theorem 1.1 in \cite{21liu2015mobius}, where the multiple ergodic means (some times called non-conventional ergodic averages)
are not concerned.
\begin{cor}\label{thm:3}
Let $X$ be a connected compact abelian group and $T$ be a minimal affine transformation
on $X$. Then $(X, T)$ has quasi-discrete spectrum of some order $d_T$ and the conclusion
\eqref{eq:1-1} in Theorem~\ref{thm:1} holds where

\indent \textnormal{(1)} $q_1, \cdots, q_\ell \in Q[x]$ are $ \ell $ $(\ge 1)$ polynomials taking values in 
$\mathbb{N}$ with $\Delta := \max_{1\le j\le \ell} d^\circ q_j$;

\indent \textnormal{(2)}  $(c_n)$ is an oscillating sequence of order $d_T\Delta$ such that $\sum_{n=1}^{N}|c_n|=O(N)$.
\end{cor}

The polynomial dynamics on the ring $\mathbb{Z}_p$ of $p$-adic integers have been well investigated
(see \cite{2anashin2009applied,11fan2011minimal}). Basing on the minimal decomposition theorem in 
\cite{11fan2011minimal} and on Theorem~\ref{thm:1} we can prove the following result.
\begin{cor}\label{thm:4}
Assume $X = \mathbb{Z}_p$ and $T$ is a polynomial in $\mathbb{Z}_p[z]$. Then the conclusion \eqref{eq:1-1} in
Theorem~\ref{thm:1} holds where

\indent \textnormal{(1)} $q_1, \cdots, q_l \in Q[x]$ are $ l $ $(\ge 1)$ polynomials taking values in $\mathbb{N}$ 
with $\Delta := \max_{1\le j\le l} d^\circ q_j$;

\indent \textnormal{(2)}  $(c_n)$ is an oscillating sequence of order $\Delta$ such that $\sum_{n=1}^{N}|c_n| = O(N)$.
\end{cor}
    
The same conclusion holds for dynamics defined by rational functions of good reduction
for which a minimal decomposition is established in \cite{8AHFappear}.
\medskip

Recall that Sarnak's conjecture concerns the orthogonality between the M\"{o}bius function $\mu$ and a topological
system $(X,T)$ with zero entropy:
\begin{equation}\label{eq:4-1}
\forall f \in C(X),\ \forall x \in X, \lim_{N\to\infty}\dfrac{1}{N}\sum_{n=0}^{N-1}\mu(n)f(T^n x)=0.
\end{equation}
Sarnak's conjecture are proved for some special cases, for example, nilsequences \cite{GT}, horocycle flows \cite{BSZ}, the Thue-Morse sequence \cite{MR} etc.
Using our Theorem~\ref{thm:1} 
 and the proof of Theorem 4.1   in \cite{7downarowicz2015odometers}, we can get  the following result.
\begin{cor}\label{thm:5}
	Let $(Y, S )$ be a unique ergodic topological system with quasi-discrete spectrum,
	having $\tau $ as the invariant measure. Let $(X, T)$ be a unique ergodic topological system having
	$ \sigma $ as the invariant measure. Suppose there exists a topological factor $\pi\colon X \to Y $, 
	which is also an isomorphism between the two measure-theoretic dynamical systems $(X, T,\sigma)$ and
	$(Y, S, \tau)$. Then Sarnak's conjecture holds for $(X, T)$, even  the M\"{o}bius function in~\eqref{eq:4-1}
	can be replaced 
by any fully oscillating sequence or  oscillating sequence of finite order (depending on
the order of the quasi-discrete spectrum of $(Y,S)$).
\end{cor}


	Here are some example of oscillating sequences:
	\begin{itemize}
		\item Let $ \alpha $  be an irrational number and $ d\ge 1 $ . The sequence $ (e^{2\pi i n^{d+1}\alpha})_{n\ge 0} $ is 
		oscillating of exact order $ d $, because $ n^{d+1}\alpha +P(n) $  is uniformly distributed modulo  for any
		$ P\in \mathbb{R}_d[x] $.
		\item The M\"{o}bius function $ (\mu(n))_{n\ge 0}  $ is fully oscillating 
		\cite{5davenport1937some,13hua2009additive} (see also \cite{22zhan1996exponential}).
		\item The random sequences of $ -1 $ or $ 1 $ (with respect to the symmetric Bernoulli probability)
		are fully oscillating almost surely. This is a special case of subnormal sequences.
	\end{itemize}

    Let us recall the notion of subnormality from probability theory which was introduced by
    Kahane \cite{19kahane1985some}. A real random variable $ \xi $ is said to be subnormal if
    \[
    \mathbb{E}e^{\lambda\xi}\le e^{\lambda^2/2},\quad \forall\,\lambda\in \mathbb{R}.
    \]
    By \textit{subnormal sequence} we mean a sequence of independent subnormal variables. An example
    of subnormal sequences is an i.i.d sequence of symmetric Bernoulli variables.
    
      \begin{thm}\label{thm:2}
    	A subnormal sequence is fully oscillating almost surely.
    \end{thm}



\section{Proofs 
             }
\begin{proof}[\quad Proof of Theorem 1] 
	First observe that the Banach algebra $C(X^\ell)$ of continuous functions defined on the product space 
	$X^\ell$ possesses  a dense subalgebra $\mathscr{L}_1$ where $\mathscr{L}_1$  is the subspace of linear 
	combinations of $g_1 \otimes \cdots \otimes g_\ell $ with $ g_1, \cdots, g_\ell \in C(X)$, called tensor 
	product of $ g_1,\cdots, g_\ell $, which are defined by
	\[
	g_1 \otimes \cdots \otimes g_\ell (x_1,\cdots,x_\ell) = g_1(x_1)g_2(x_2)\cdots g_\ell(x_\ell).
	\]
	The density is a consequence of Stone-Weierstrass theorem. Since  quasi-eigenfunctions
	are dense in $C(X)$ by (ii), the algebra $C(X^\ell)$ possesses another dense subspace $\mathscr{L}_2$, which consists
	of linear combinations of tensor products of quasi-eigenfunctions. Notice that $ \mathscr{L}_2$ is actually
	a subalgebra, because the quasi-eigenfunctions form a group.
	
	On the other hand, since $\sum_{n}^{N-1}|c_n| = O(N) $, there exists a constant $ D > 0 $ such that
	\[
	\forall F\in C(X^\ell),\qquad \left\|\dfrac{1}{N}\sum_{n=0}^{N-1}c_n F(T^{q_1(n)},\cdots,T^{q_\ell(n)})
	\right\|_\infty \le D\|F\|_\infty.
	\]
	So, using an approximation argument, we have only to prove~\eqref{eq:1-1} for $ F = g_1 \otimes \cdots \otimes g_\ell$,
	where $ g_1,\cdots,g_\ell $ are quasi-eigenfunctions.
	
	Let $ f $ be a quasi-eigenfunction of order $ k $ ($\le d$). Then there exist $ f_0 \in G_0
	$, $f_1 \in G_1,\cdots,  f_{k-1} \in G_{k-1}, f_k = f $ such that
	\[
	f_k \circ T = f_{k-1}f_k,\ \cdots ,\  f_2\circ T = f_1 f_2,\  f_1 \circ T = f_0  f_1.
	\]
	Then, for each integer $n \ge 0$, we have
	\begin{equation}\label{eq:2-1} 
	f (T^n x) = f_0(x)^{\binom{n}{k}}f_1(x)^{\binom{n}{k-1}}\cdots f_{k-1}(x)^{\binom{n}{1}} f_k(x)^{\binom{n}{0}}.
	\end{equation}
	(Remark that $ \binom{n}{k}= 0 $ for $ 0 \le n < k$ ). This can be easily proved by induction on $n$ using
	Pascal's formula. Fix $ x \in X $. Let $ \theta_j \in [0, 1) $ be the number such that $ f_j(x) = e^{2\pi i 
		\theta_j},0 \le j \le k $. Then we can write the equality \eqref{eq:2-1} as
	\begin{equation}\label{eq:2-2}
	\forall n \ge 0,\quad f(T^n x) = e^{2\pi iQ(n)}
	\end{equation}
	where
	\[
	Q(z) = \sum_{j=0}^{k}\theta_j\binom{z}{k-j} \in \mathbb{R}_k[z].
	\]
	Now assume $ F = g_1 \otimes \cdots \otimes g_\ell $, where $ g_j $ is a quasi-eigenfunction of order $k_j$ 
	($1 \le j \le \ell$). Applying~\eqref{eq:2-2} to each $g_j$, we get a polynomial $Q_j \in \mathbb{R}_{k_j}[z] $
	such that
	\[
	g_j(T^nx) = e^{2\pi iQ_j(n)}.
	\]
	Therefore
	\[
	F(T^{q_1(n)} x,\cdots,T^{q_\ell(n)} x) = g_1(T^{q_1(n)}x) \cdots g_\ell(T^{q_\ell(n)}x) = e^{2\pi iP(n)}
    \]
    where $$
     P = \sum_{i=1}^{\ell}Q_j\circ q_j \in \mathbb{R}_{d\Delta}[z] $$
     (recall that $ \Delta = \max(d^\circ q_1, 
     \cdots, d^\circ q_\ell)$). Thus
    \[
     \dfrac{1}{N}\sum_{n=0}^{N-1}c_n F(T^{q_1(n)},\cdots,T^{q_\ell(n)}) = 
     \dfrac{1}{N}\sum_{n=0}^{N-1}c_n e^{2\pi i P(n)}
    \]
which tends to zero, by (iii).
\end{proof}

If $(X, T)$ is not of quasi-discrete spectrum, the tensor products of quasi-eigenfunctions still
span a subalgebra of $ C(X^\ell) $. Its closure is a closed subspace of $ C(X^\ell) $, that we denote 
by $C_T(X^\ell)$. The conclusion of Theorem~\ref{thm:1} remains true for all $F$ in $C_T(X^\ell)$.

\begin{proof}[\quad Proof of Corollary 2] 
By the main theorem in \cite{11fan2011minimal}, $ \mathbb{Z}_p = M \cup B $ (disjoint union) 
where $ M $ is a disjoint union of at most countably many minimal sets $M_i$ and every point in $B$ is mean
attracted to one of minimal sets in $M$ (see \cite{10fan2015oscillating} for the definition of mean attraction), and
each subsystem $ T\colon M_i \to M_i $ is either a finite cycle or is conjugate to an adding machine.
So, $ T \colon M_i \to M_i $ has discrete spectrum, i.e. quasi-discrete spectrum of order $1$. Thus 
\eqref{eq:1-1} holds for all $ x \in M_i$, by Theorem~\ref{thm:1}. For $ x \in B $, we can prove \eqref{eq:1-1} 
as in \cite{10fan2015oscillating} (see the end of the proof of Theorem 1 in \cite{10fan2015oscillating}).
\end{proof}

The same result in Corollary~\ref{thm:4} holds for $p$-adic dynamics defined by rational functions with
good reduction, because a minimal decomposition is also proved for such dynamics in [9].
Recall that any rational function $\phi \in \mathbb{Q}_p[z]$ having good reduction is $1$-Lipschitz continuous
on the projective line $\mathbb{P}^1(\mathbb{Q}_p)$:
\[
\rho(\phi(P_1),\phi(P_2)) \le \rho(P_1, P_2)\quad \text{for all $P1, P2 \in \mathbb{P}^1(\mathbb{Q}_p)$}
\]
where $\rho(\cdot, \cdot)$ is the \textit{spherical metric} on $\mathbb{P}^1(\mathbb{Q}_p)$ 
(see \cite[p.59]{26silverman2007arithmetic}).	

\begin{proof}[\quad Proof of Theorem~\ref{thm:2}]
  Consider a subnormal sequence $ (\xi_n) $ on a probability space $ (\Omega, \mathscr{A}, \mathbb{P}) $ and a finite 
    sequence of trigonometric polynomials of $ s $ variables $ f_n(t_1,\cdots,t_s) $ of order less than or equal to N.
    Then consider the random trigonometric polynomial
    \[
    P(t) = \sum_{n}\xi_n f_n(t).
    \]
    There exists an absolute constant $C > 0$ such that
    \begin{equation}\label{eq:ineq1-3}
    \mathbb{P}\left(\|P\|_\infty \le C\sqrt{s\sum\|f_n\|_\infty^2\log N} \right) \le \dfrac{1}{N^2}.
    \end{equation}
    We refer this inequality as Littlewood-Salem-Kahane inequality (\cite{19kahane1985some} p. 70, Theorem
    3. See \cite{12fan2003inegalite} for a different proof and a generalization). Basing on this inequality, we
    can prove the conclusion as follows.

	We have to prove that for almost all $ \omega \in \Omega $  and for all polynomials $ P \in \mathbb{R}[z]$, 
	we have
	\begin{equation}\label{eq:3-1} 
	\lim_{N\to\infty} \dfrac{1}{N} \sum_{n=0}^{N-1}\xi_n(\omega) e^{2\pi iP(n)} = 0.
	\end{equation}
	Let $d = d^\circ P$ be the degree of $ P $ and $ P(x) = \sum_{j=0}^d t_jx^j $ with $(t_0, t_1,\cdots, t_d) \in 
	\mathbb{R}^{d+1} $. Then
	\[
	e^{2\pi iP(n)} = e^{2\pi i(t_0+ t_1 n+ \cdots + t_dn^d)} =: f_n(t_0,t_1,\cdots,t_d).
	\]
	We regard $f_n$ as a trigonometric polynomial of $d + 1$ variables, of degree $1 + n + \cdots + n^d \le 
	(d + 1)N^d$. Thus 
	\[
		\sup\limits_{P\in \mathbb{R}_d[z]}\left|\sum_{n=0}^{N-1}\xi_n(\omega)e^{2\pi P(n)}\right| = 
		\max\limits_{t\in \mathbb{T}^{d+1}} \left|\sum_{n=0}^{N-1}\xi_n(\omega)f_n(t)\right|.
	\]
	Observe that $ \|f_n\|_{\infty} = 1 $. By Littlewood-Salem-Kahane inequality $\eqref{eq:ineq1-3}$, we have
	\[
	\mathbb{P}\left( \sup\limits_{P\in \mathbb{R}_d[z]}\left|\sum_{n=0}^{N-1}\xi_n(\omega)e^{2\pi P(n)}\right| 
	\ge C(d+1)\sqrt{(d+1)N\log N}\right)\le \dfrac{1}{N^2}.
	\]
	Then, by Borel-Cantelli lemma, for almost all $\omega$,
	\[
	\sup\limits_{P\in \mathbb{R}_d[z]}\left|\sum_{n=0}^{N-1}\xi_n(\omega)e^{2\pi P(n)}\right| = O(\sqrt{N\log N}).
	\]
	Since $d$’s are countable, this implies \eqref{eq:3-1}. 
\end{proof}

We have actually proved that almost surely
\[
\dfrac{1}{N}\sum_{n=0}^{N-1}\xi_n(\omega)e^{2\pi i P(n)} = O_{d,\omega}\left(\sqrt{\dfrac{\log N}{N}} \right),
\]
where the constant involved in $O$ is uniform for all polynomials of degree less than some
fixed $d$. In other words, it depends only on the degree $d$ and on $\omega$.

\section{Some remarks}

1. Hahn and Parry \cite{14hahn1965minimal} restricted  their study of quasi-discrete systems to minimal
homeomorphisms and they made the assumptions that the system is  totally minimal and  has no non-trivial
eigenvalues of finite order.  For our concern, it is not necessary to make these assumptions.

Assume the space $X$ is connected. Then the minimality of $T$ implies
automatically the total minimality. 
Under the
extra assumption that there is no non-trivial eigenvalues of finite order, a minimal homeomorphism
of quasi-discrete spectrum must be uniquely ergodic and totally ergodic (Theorem~1 and Theorem~2 
in \cite{14hahn1965minimal}), and two minimal homeomorphisms of quasi-discrete spectrum are conjugate iff 
their systems of quasi-eigenvalues are equivalent (Theorem 3 in \cite{14hahn1965minimal}). It is easy to see 
that the total minimality implies no eigenvalues of finite order. Then $(X,T)$ has no non-trivial eigenvalues 
of finite order if $ X $ is connected and $T$ is minimal. 

We don't assume that $X$ is connected. This allows us 
to include interesting dynamics on disconnected spaces. Here is an illustrative example. Let us consider
an affine system  $ Tx = ax+b $ on the ring $\mathbb{Z}_p$ of $p$-adic integers, with $ a,b \in \mathbb{Z}_p $. 
It is well known that $T$ is minimal iff $a = 1 \mod p$ and $ b \ne 0 \mod p $ 
(see \cite{2anashin2009applied,9fan2007strict}). Under the assumption $ a = 1 \mod p $ and $ b \ne 0 \mod p$, 
the dynamics $T$ is minimal and uniquely ergodic and it has discrete spectrum [10]. But $T$ is not totally 
minimal and not totally ergodic either (because $T^m$ is not minimal when $ m = 0 \mod p $), and all the eigenvalues
 of $T$ are of finite order! \vspace{0.5em}

2. It was also proved by Hahn and Parry \cite{14hahn1965minimal} that under the hypothesis that there is no
non-trivial eigenvalues of finite order, a minimal homeomorphism of quasi-discrete spectrum
can be represented as an affine transformation of a compact abelian group. Theorem~2 in \cite{16hoare1966affine}
states that if there exists a totally minimal affine transformation on a compact abelian group,
the group must be connected. Hoare and Parry  proved in \cite{17hoare1966affine2} that all minimal affine
transformations of any connected compact abelian group have quasi-discrete spectrum, to which Theorem 1 applies (Corollary 1). \vspace{0.5em}

3. Dynamics on disconnected spaces like the symbolic space $\{0, 1,\cdots,m-1\}$ ($m \ge 2$) are worthy of study. 
But such dynamics were excluded in the works \cite{14hahn1965minimal,16hoare1966affine,17hoare1966affine2}. 
The dynamics of polynomial dynamics on $\mathbb{Z}_p$ have been well investigated (see 
\cite{2anashin2009applied,11fan2011minimal}). 
If such a polynomial is minimal, then it is conjugate to the adding machine $ x \mapsto x + 1 $
and therefore it has not only quasi-discrete spectrum but discrete spectrum and the eigenvalues are all of
finite order. Polynomials in $\mathbb{Z}_p[z]$ and rational functions in $\mathbb{Q}_p[z]$ of good reductions 
share a minimal decomposition \cite{11fan2011minimal,8AHFappear}. Therefore they locally behave like an 
adding machine (in general adding machine on $\prod_{k=1}^\infty \mathbb{Z}/m_k\mathbb{Z}$ for some sequence of 
integers $(m_k)$).\vspace{0.5em}

4. Under the assumption that $(X, T)$ is a uniquely ergodic model of some totally ergodic
measure-theoretic system $(Y, \sigma, S )$ with quasi-discrete spectrum in the sense of Abramov [1],
Abdalaoui, Lema\'{n}czyk and de la Rue proved (Theorem 1 in [3]) that for any $ f \in C(X)$ and
any $ x \in X$ we have
\[
\lim\limits_{N\to \infty}\dfrac{1}{N}\sum_{n=0}^{N-1}\nu_{n}f(T^n x) = 0
\]
where $\nu\colon \mathbb{N} \to C $ is a {\em multiplicative function} with $|\nu(n)| \le 1$ and $\sum_{n=0}^{N}\nu(n)=
o(N)$. Actually the condition of quasi-discrete spectrum in the sense of Abramov can be replaced by the
weaker condition of Asymptotic Orthogonal Powers (Theorem 2 in \cite{3lemanczyk2015}): for any $f, g \in 
L^2(\sigma)$ with $\int f\,d\sigma = \int g\,d\sigma=0 $, we have
\[
\lim_{p,q\to\infty}\sup_{\kappa\in J^e(S^p,S^q)}\left|\int_{Y\times Y}f\otimes g\,d\kappa\right| = 0
\]
where $p, q$ are different primes and $J^e(S^p,S^q)$ is the set of ergodic joinings of $S^p$ and $S^q$. 
In this setting, the total ergodicity plays an important role in the arguments by joining and the
proof is based on Bourgain-Sarnak-Ziegler orthogonality criterion \cite{BSZ}.
Theorem 1 in [3] and Corollary \ref{thm:5} in the present note have some overlaps, but no one can implies the other.

\vspace{0.5em}

5. Huang,Wang and Zhang \cite{18huang2016m} proved that M\"{o}bius sequence is orthogonal to all topological
models of an ergodic system with discrete spectrum, answering affirmatively a question
raised by Downarowicz and Glasner \cite{6downarowicz2015isomorphic}. The result recovers and generalizes some results
obtained in \cite{10fan2015oscillating} and \cite{27Wei2016entropy}. Actually they found a sufficient condition for one 
point observable to be orthogonal to the M\"{o}bius sequence. The proof uses the estimate on short averages of
non-pretentious multiplicative functions due to Matom\"{a}ki, Radziwill and Tao \cite{23matomaki2015averaged}.\vspace{0.5em}

\textit{Acknowledgment.}  I would like to thank the support of the Wallenberg Foundation and
the Royal Swedish Academy of Sciences, and the hospitality of the Centre for Mathematical
Sciences in University of Lund.   
\bibliographystyle{plain}

\begin{thebibliography}{9}	
\bibitem{1abramov1962metric} L. M. Abramov,\,\emph{Metric automorphisms with quasi-discrete spectrum}, 
Izv. Akad. Nauk. U. S. S. R., 26(1962),513-530.

\bibitem{2anashin2009applied} V. Anashin and A. Khrennikov,\,\emph{Applied algebraic dynamics}, 
de Gruyter Expositions in Mathematics 49, Walter de Gruyter, 2009.

\bibitem{3lemanczyk2015} H. El Abdalaoui, M. Lema\'{n}czyk and T. de la Rue,\,\emph{Automorphisms with quasi-discrete spectrum,
multiplicative functions and average orthogonality along short intervals}.
\url{https://hal.archives-ouvertes.fr/hal-01176039/document}, 2015.

 \bibitem{B} J. Bourgain, \, \emph{ On the correlation of the Moebius function with rank-one systems}, J. Anal. Math. 120 (2013), 105-130.

\bibitem{BSZ} J. Bourgain, P. Sarnak, and Ziegler,\,\emph{Disjointness of M\"{o}bius from horocycle flows. From Fourier 
analysis and number theory to Radon transforms and geometry}, 67-83, Dev. Math. 28, Springer, New York, 2013.

\bibitem{5davenport1937some} H. Davenport,\,\emph{On some infinite series involving arithmetical 
	functions (II)}, Quart. J. Math. Oxford, 8(1937), 313-320.

\bibitem{6downarowicz2015isomorphic} T. Downarowicz and E. Glasner,\,\emph{Isomorphic extensions and applications}, 
Topological Methods in Nonlinear Analysis (2015), to appear, arXiv:1502.06999v1.

\bibitem{7downarowicz2015odometers} T. Downarowicz and S. Kasjan,\,\emph{Odometers and Toeplitz systems revisited 
	in the context of Sarnak's conjecture}, 2015. arXiv: 1502.02307v1		

\bibitem{8AHFappear} A. H. Fan, S. L. Fan, L. M. Liao and Y. F. Wang,\,\emph{On minimal decomposition of p-adic rational 
functions with good reduction}, Discrete and Continuous Dynamical Systems, to appear.

\bibitem{9fan2007strict} A. H. Fan, M. T. Li and J. Y. Yao D. Zhou,\,\emph{Strict ergodicity of affine $ p $-adic 
	dynamical systems on $ Z_p $}, Adv. Math., 214 (2007) 666-700.

\bibitem{10fan2015oscillating} A. H. Fan and Y. P. Jiang,\,\emph{Oscillating sequences, minimal mean attractability 
	and minimal mean-Lyapunov-stability}. Erg. Th. Dynam. Syst., to appear.

\bibitem{11fan2011minimal} A. H. Fan and L. M. Liao,\,\emph{On minimal decomposition of p-adic polynomial dynamical 
	systems}, Adv. Math., 228 (2011), 2116-2144.

\bibitem{12fan2003inegalite} A. H. Fan and D. Schneider,\ \emph{Sur une in\'{e}galit\'{e} de Littlewood-Salem}, 
Ann. I. H. Poincar\'{e} PR 39, 2 (2003),193-216.

\bibitem{GT} B. Green and T. Tao, \, \emph{The M\"{o}bius function is strongly orthogonal to nilsequences },  Ann. of Math. (2) 175, no 2 (2012), 541-566.

\bibitem{13hua2009additive} L. G. Hua,\,\emph{Additive Theory of Prime Numbers} (Translations of Mathematical Monographs 
: Vol 13).Amer Mathematical Society. 1966.

\bibitem{14hahn1965minimal} F. Hahn and W. Parry,\,\emph{Minimal dynamical systems with quasi-discrete spectrum}. 
J. London Math. Soc., \textbf{40}\,(1965), 309-323.

\bibitem{15halmos1942operator} P. R. Halmos and J. von Neumann,\,\emph{Operator methods in classical mechanics. II}, 
Ann. of Math. (2) 43\,(1942), 332-350.

\bibitem{16hoare1966affine} H. Hoare and W. Parry,\,\emph{Affine transformations with quasi-discrete spectrum (I)}. 
J. London Math. Soc., \textbf{41}\,(1966), 88-96.

\bibitem{17hoare1966affine2} H. Hoare and W. Parry,\,\emph{Affine transformations with quasi-discrete spectrum (II)}. 
J. London Math. Soc., \textbf{41}\,(1966), 529-530.

\bibitem{18huang2016m} W. Huang, Z. R. Wang and G. H. Zhang,\,\emph{M\"{o}bius disjointness for topological models of 
	ergodic systems with discrete spectrum}, arXiv:1608.08289v2

\bibitem{19kahane1985some} J. P. Kahane,\,\emph{Some random series of functions}, Cambridge University Press, 1985.

\bibitem{20li2015mean} J. Li, S. Tu and X. D. Ye,\,\emph{Mean-equicontinuity and mean sensitivity},
Ergodic Theory Dynam. Systems 35\,(2015), no. 8, 2587-2612.

\bibitem{21liu2015mobius} J. Y. Liu and P. Sarnak,\,\emph{The M\"{o}bius function and distal flows}, 
Duke Math. J. 164(2015),no.\,7,1353-1399.

\bibitem{22zhan1996exponential} J. Y. Liu and T. Zhan,\,\emph{Exponential sums involving the M\"{o}bius function}. 
Indagationes Mathematicae Volume 7 (1996), Issue 2, 271-278.

\bibitem{23matomaki2015averaged} K. Matom\"{a}ki, M. Radziwill, and T. Tao,\,\emph{An averaged form of Chowla's 
	conjecture}, Algebra Number Theory, 9 (2015), no. 9, 2167-2196.
	
\bibitem{MR}	C.Mauduit and J.Rivat, \, \emph{Sur un probl\`eme de Gelfond: la somme des chiffres des nombres premiers}, Ann.
of Math. (2) 171 (2010), 1591-1646.

\bibitem{24Sarnak2009three} P. Sarnak,\,\emph{Three lectures on the M\"{o}bius function, randomness and dynamics}, 
IAS Lecture Notes, 2009;\url{ http://publications.ias.edu/sites/default/files/MobiusFunctionsLectures(2).pdf}.

\bibitem{25sarnak2010mobius} P. Sarnak,\,\emph{M\"{o}bius randomness and dynamics}, 
Not. S. Afr. Math. Soc. 43 (2012), 89-97.

\bibitem{26silverman2007arithmetic} J. Silverman,\,\emph{The arithmetic of dynamical systems}, volume 241 of Graduate 
Texts in Mathematics. Springer, New York, 2007.

\bibitem{27Wei2016entropy} F.\,Wei,\,\emph{Entropy of arithmetic functions and Sarnak's M\"{o}bius disjointness 
	conjecture}, 2016. Thesis (Ph.D.), The University of Chinese Academy of Sciences.
\end{thebibliography}

\end{document}